\documentclass[11pt]{amsart}
\usepackage[dvipdfmx]{graphicx}
\usepackage{graphicx}
\usepackage{amssymb,amstext,amsfonts,amscd}
\usepackage{amsmath}
\usepackage{latexsym}
\usepackage{amsthm}
\usepackage{multicol}
\usepackage{fancybox}
\usepackage{here}
\usepackage{multirow}
\usepackage{mathrsfs}

\usepackage{enumitem}

\makeatletter
\@namedef{subjclassname@2020}{%
  \textup{2020} Mathematics Subject Classification}
\makeatother

\def\opn#1#2{\def#1{\operatorname{#2}}}
\opn\max{max}
\opn\min{min}
\opn\rank{rank}
\opn\Ker{Ker}
\opn\id{id}
\opn\mod{mod}
\opn\ord{ord}
\opn\det{det}
\opn\Cone{Cone}
\opn\Int{Int}
\opn\grad{grad}
\opn\codim{codim}
\opn\sign{sign}

\newcommand{\bz}{\mathbb{Z}}

\newcommand{\br}{\mathbb{R}}
\newcommand{\bc}{\mathbb{C}}

\newcommand{\zero}{\boldsymbol{0}}
\newcommand{\ba}{\boldsymbol{a}}   %%%% boldsymbol

\newcommand{\bw}{\boldsymbol{w}}
\newcommand{\bx}{\boldsymbol{x}}

\newcommand{\bbz}{\boldsymbol{z}}

\newcommand{\QED}{$\Box$}

\newtheorem{theorem}{Theorem}
\newtheorem{definition}[theorem]{Definition}

\newtheorem{proposition}[theorem]{Proposition}
\newtheorem{lemma}[theorem]{Lemma}

\newtheorem{example}[theorem]{Example}

\numberwithin{equation}{section}

\newcommand{\bba}{\boldsymbol{\alpha}}
\newcommand\rdeg{{\rm{rdeg}\/}}
\newcommand\pdeg{{\rm{pdeg}\/}}
\usepackage{multirow}

\begin{document}
\pagestyle{plain}

\title
{A note on Newton non-degeneracy of mixed weighted homogeneous polynomials}

\author{Sachiko Saito and Kosei Takashimizu}

\address{
Department of Mathematics Education, Asahikawa Campus, Hokkaido University of Education, 
Asahikawa 070-8621, Hokkaido, Japan,\ \ 
Seiryo Junior High School, 
Kushiro 085-0814, Hokkaido, Japan. %% 釧路市立青陵中学校教諭
}

\email{saito.sachiko@a.hokkyodai.ac.jp,\ tkosei25@gmail.com}

% 14P05 Real algebraic sets
% 32S45 Modifications; resolution of singularities (complex-analytic aspects)
% 32S55 Milnor fibration; relations with knot theory


\begin{abstract}
A mixed polynomial $f(\bbz, \bar{\bbz})$ is called a mixed weighted homogeneous polynomial (Definition \ref{mixed-whp}) 
if it is both radially and polar weighted homogeneous. 
%%%
Let $f$ be a mixed weighted homogeneous polynomial with respect to 
a strictly positive radial weight vector $P$ and a polar weight vector $Q$. 
Suppose that $f$ is Newton non-degenerate 
over a compact face $\Delta(P)$ %%% 
and polar weighted homogeneous of non-zero polar degree with respect to $Q$. 
Then $f : {\bc^*}^n \to \bc$ has no mixed critical points. 
Moreover, under the assumption $f^{-1}(0) \cap {\bc^*}^n \neq \emptyset$, 
$f : {\bc^*}^n \to \bc$ is surjective. 
In other words, in this situation, 
Newton non-degeneracy over a compact face $\Delta(P)$ implies 
strong Newton non-degeneracy over $\Delta(P)$ (Proposition \ref{remark4}). 
With this fact as a starting point, we investigate the sets $f^{-1}(0) \cap {\bc^*}^n$, and 
show the existence of 
a collection of mixed weighted homogeneous polynomials $f = f_{\Delta (P)}$ of non-zero polar degree 
which satisfy $\dim \Delta (P) \geq 1$ and $f^{-1}(0) \cap {\bc^*}^n = \emptyset$ (Theorem \ref{true-nondegenerate-Thm}). 
We also give an example of 
convenient mixed function germs 
of mixed weighted homogeneous face type which are not true non-degenerate (Definition \ref{true-non-degenerate}). 
\end{abstract}

\maketitle

%%%%%%%%%%%%%%%%%%%%%%%%%%%%%%%%%%%%%%%%%%%%%%%%%%%%

\section{Mixed functions and their radial Newton polyhedrons}\label{section-mixed-function}

Let $U$ be a neighborhood of $\zero$ in $\bc^n$. 
We assume $\bar{U}=U$, 
where $\bar{\bbz}$ stands for 
the complex conjugate $(\bar{z_1}, \dots , \bar{z_n})$ of $\bbz=(z_1, \dots , z_n) \in \bc^n$. 
For a complex valued holomorphic function $F(\bbz, \bw)$ on $U\times U$ with complex $2n$ variables, 
we set 
$$
f(\bbz, \bar{\bbz}) := F(\bbz, \bar{\bbz}), 
$$
which is defined over $U$. 
We call it 
a {\em mixed analytic function} (or {\em mixed function}) $f$ on $U$. 
We assume that 
$
F(\zero, \zero) = 0
$
throughout this paper, 
and hence, 
$f(\zero)=0$. 
Let 
$
F(\bbz, \bw) = \sum_{\nu, \mu}c_{\nu, \mu} \bbz^\nu \bw^\mu
$
be the Taylor expansion of $F$ at $(\zero , \zero)$, 
where 
$
\nu = (\nu_1, \dots , \nu_n),\ \mu = (\mu_1, \dots , \mu_n),\ 
\nu_i \geq 0,\ \mu_j \geq 0,\ \ 
\bbz^\nu := z_1^{\nu_1}\cdots z_n^{\nu_n},\ \ \bw^\mu := w_1^{\mu_1}\cdots w_n^{\mu_n} .
$
Then we have 
\begin{equation}\label{Taylor}
f(\bbz, \bar{\bbz}) = \sum_{\nu, \mu}c_{\nu, \mu} \bbz^\nu \bar{\bbz}^\mu .
\end{equation}
Note that the coefficients of the Taylor expansion \eqref{Taylor} of a mixed function $f$ are unique. %%%
We call $f(\bbz, \bar{\bbz})$ a {\em mixed polynomial} 
if 
the number of monomials $c_{\nu, \mu} \bbz^\nu \bar{\bbz}^\mu,\ c_{\nu, \mu}\neq 0$ is finite. 

Let $f(\bbz, \bar{\bbz})$ be a mixed function on $U \ (\subset \bc^n)$. 
We say 
$
\ba = (a_1, \dots , a_n) \in U
$
is a {\em mixed critical point} (or a {\em mixed singular point}) of $f$ 
if 
the rank of the differential map 
$(df)_{\ba} : T_{\ba} \bc^n \to T_{f(\ba)} \bc \cong T_{f(\ba)} \br^2$ 
is less than $2$. 
We say $\ba \in U$ is a {\em mixed regular point} of $f$ 
if it is not a mixed critical point of $f$. 

We set 
$K_+^n := \{ (x_1, \dots , x_n) \in K^n \ |\ x_i \geq 0\}$, where $K=\br \ \text{or} \ \bz$. 
Let $(f,\zero)$ be the germ of a mixed function 
$
f(\bbz, \bar{\bbz}) = \displaystyle \sum_{\nu, \mu}c_{\nu, \mu} \bbz^\nu \bar{\bbz}^\mu
$
at $\zero \in \bc^n$. 
Let 
$
\Gamma_+(f)
$
be the convex hull of the set 
$
\displaystyle \bigcup_{c_{\nu, \mu} \neq 0}  (\nu + \mu) + \br_+^n ,
$
which we call 
the {\em (radial) Newton polyhedron} of the germ $(f,\zero)$ of a mixed function $f$ at $\zero$. 

For a ``weight vector" %% 整数ベクトル
$
P={}^t(p_1, \dots , p_n)\ (\neq \zero) \in \bz_+^n
$, 
let 
$d(P)$ 
be the minimum value of the linear function 
$
P : \Gamma_+(f) \to \br,\ P(\xi) := \displaystyle \sum_{j=1}^n p_j \xi_j,
$
where 
$\xi = (\xi_1 , \dots , \xi_n) \in \Gamma_+(f)$. 
We set 
$
\Delta(P) := \{ \xi \in \Gamma_+(f) \ |\ P(\xi) = d(P) \}, 
$
which we call a {\em face} of $\Gamma_+(f)$. 
Note that $\Delta(P) \neq \emptyset$ by its definition. 
%%% If $\dim \Delta(P) = n-1$, then $P$ is unique up to multiplications of positive real numbers. 

We say a weight vector $P={}^t(p_1, \dots , p_n) \in \bz_+^n$ is {\em strictly positive} if 
$p_i >0$ for every $i \ (=1, \dots , n)$. 
We write $P\gg 0$ if $P$ is strictly positive. 
Note that a face $\Delta$ of $\Gamma_+(f)$ is compact 
if and only if 
$\Delta = \Delta(P)$ 
for some strictly positive weight vector $P$ (see Lemma \ref{compact-face} in \S \ref{section-appendix}). 

For a \underline{compact} face $\Delta(P)$, we define 
$$
f_{\Delta(P)}(\bbz) \ (\text{or}\ f_P(\bbz)\,) := \sum_{\mu + \nu \in \Delta(P)}c_{\nu, \mu} \bbz^\nu \bar{\bbz}^\mu ,
$$
which we call a {\em face function {\rm (or}\ face polynomial{\rm )}} 
of a mixed function germ $(f,\zero)$ (\cite{Oka2018}, p.78). 

\medskip

\section{Radially and polar weighted homogeneous polynomials}\label{section-whp}

\begin{definition}[\cite{Oka2015};\ \cite{Oka2018}, p.182]
{\rm \ \ 

\begin{itemize}
\item A mixed polynomial 
$f(\bbz, \bar{\bbz}) = \sum_{\nu, \mu}c_{\nu, \mu} \bbz^\nu \bar{\bbz}^\mu$ is called 
{\em radially weighted homogeneous} if there exists a weight vector 
$P={}^t(p_1, \dots , p_n)\ (\neq \zero) \in \bz_+^n$   %%% 非負整数ベクトル
and 
a positive integer $d_r \ (> 0)$ such that 
$
c_{\nu, \mu} \neq 0 \Longrightarrow P(\nu + \mu) = \sum_{i=1}^n p_i(\nu_i + \mu_i) = d_r .
$
We call $d_r$ {\em the radial degree} of $f$, and define 
$
\rdeg_P f := d_r .
$

\item A mixed polynomial $f(\bbz, \bar{\bbz}) = \sum_{\nu, \mu}c_{\nu, \mu} \bbz^\nu \bar{\bbz}^\mu$ is called 
{\em polar weighted homogeneous} if there exists a weight vector 
$Q={}^t(q_1, \dots , q_n)\ (\neq \zero) \in \bz^n$ 
%%%% \cite{Oka2018}p.183, 注意8.6では，p_iが負である場合も考えている
%%%%% 後で用いるのは「極次数」が0でないということ
and 
an integer $d_p$ ($>0,\ 0\ \text{or}\ <0$) such that 
$
c_{\nu, \mu} \neq 0 \Longrightarrow Q(\nu - \mu) = \sum_{i=1}^n q_i(\nu_i - \mu_i) = d_p .
$
We call $d_p$ {\em the polar degree} of $f$, and define 
$
\pdeg_Q f := d_p .
$
\end{itemize}
}
\end{definition}

Since we assume $f(\zero)=0$, every face function $f_{\Delta(P)}(\bbz)$, 
where $P$ is strictly positive, of a mixed function germ $(f,\zero)$ 
is a radially weighted homogeneous polynomial of radial degree $d(P)\ (>0)$ 
with respect to the weight vector $P$. 

\begin{example}[\cite{Oka2018}, Example 9.17]
$f(\bbz, \bar{\bbz}) := z_1^2 \bar{z_1} - z_2 \bar{z_2}^2$ is radially weighted homogeneous with respect to $P={}^t(1,1)$ and 
polar weighted homogeneous with respect to $Q={}^t(1,-1)$. 
\end{example}

\medskip

For radially and polar weighted homogeneous polynomials, we have the following facts: 
\begin{lemma}[cf.\cite{Oka2008}, \cite{Oka2010}]\label{mixed-action}
Let $f(\bbz, \bar{\bbz})$ be a mixed polynomial. We have the following. 
\begin{itemize}
\item Let $P={}^t(p_1, \dots , p_n)\ (\neq \zero) \in \bz_+^n$ be a weight vector and $d_r$ be a positive integer. 
For a positive real number $t$ and $\bbz \in \bc^n$, we define 
$
t \circ \bbz := (t^{p_1} z_1, \dots , t^{p_n} z_n)
$. 
Then 
$
f(t \circ \bbz) = t^{d_r}f(\bbz)\ \ \ (\forall t>0, \forall \bbz \in \bc^n)
$
if and only if 
$
c_{\nu, \mu} \neq 0 \Longrightarrow P(\nu + \mu) = \sum_{i=1}^n p_i(\nu_i + \mu_i) = d_r
$. 
\item Let $Q={}^t(q_1, \dots , q_n)\ (\neq \zero) \in \br^n$ %%%%
be a weight vector and $d_p$ be an integer. 
For a real number $\theta$ and $\bbz \in \bc^n$, we define 
$
\theta \circ \bbz := (e^{iq_1 \theta} z_1, \dots , e^{iq_n \theta} z_n)
$. 
Then 
$
f(\theta \circ \bbz) = e^{id_p \theta}f(\bbz)\ \ \ (\forall \theta \in \br, \forall \bbz \in \bc^n)
$
if and only if 
$
c_{\nu, \mu} \neq 0 \Longrightarrow Q(\nu - \mu) = \sum_{i=1}^n q_i(\nu_i - \mu_i) = d_p
$. \hfill \QED
\end{itemize}
\end{lemma}

\begin{proposition}[Euler equalities,\ \cite{Oka2018},\ \cite{Takashimizu2021thesis},\ \cite{Saito-Takashimizu2021-8}]\label{Euler}
Let $f(\bbz,\overline{\bbz})=\sum_{\nu,\mu}c_{\nu,\mu}\bbz^{\nu}\overline{\bbz}^{\mu}$ be a mixed polynomial. 
\begin{description}
\item[(R)]\ If $f(\bbz,\overline{\bbz})$ is a radially weighted homogeneous polynomial 
of radial degree $d_{r}\ (>0)$ with respect to a weight vector $P={}^{t}\!(p_{1},\dots,p_{n})$, then we have 
\begin{equation}\label{reulereq}
\sum_{j=1}^{n}p_{j}\left(z_{j}\frac{\partial f}{\partial z_{j}} + \overline{z_{j}}\frac{\partial f}{\partial \overline{z_{j}}} \right)
= 
d_{r}f(\bbz,\overline{\bbz}) .
\end{equation}
\item[(P)]\ If $f(\bbz,\overline{\bbz})$ is a polar weighted homogeneous polynomial 
of polar degree $d_{p}$ with respect to a weight vector $Q={}^{t}\!(q_{1},\dots,q_{n})$, then we have 
\begin{equation}\label{peulereq}
\sum_{j=1}^{n}q_{j}\left(z_{j}\frac{\partial f}{\partial z_{j}} - \overline{z_{j}}\frac{\partial f}{\partial \overline{z_{j}}} \right)
= 
d_{p}f(\bbz,\overline{\bbz}) . \ \ \Box
\end{equation}
\end{description}
\end{proposition}
Although Lemma \ref{mixed-action} and Proposition \ref{Euler} are well known, 
we would like to give the proofs of them in \S \ref{section-appendix}. 

\medskip

\begin{definition}[mixed weighted homogeneous polynomial,\ Oka \cite{Oka2018}, pp.182--184]\label{mixed-whp}  %% polynomial
{\rm 
We say a mixed polynomial $f(\bbz, \bar{\bbz})$ is 
a {\em mixed weighted homogeneous polynomial}
if 
it is both radially and polar weighted homogeneous. 
Here, the corresponding weight vectors $P$ and $Q$ are possibly different. 
%%%
We say a mixed weighted homogeneous polynomial $f(\bbz, \bar{\bbz})$ is 
a {\em strongly mixed weighted homogeneous polynomial} if 
$f$ is radially and polar weighted homogeneous 
with respect to the \underline{same} weight vector $P$. 
%%%
Furthermore, 
a mixed weighted homogeneous polynomial $f$ is called 
a {\em strongly polar positive mixed weighted homogeneous polynomial} 
with respect to a weight vector $P$ 
if 
$f$ is radially and polar weighted homogeneous with respect to the same weight vector $P$ and 
$\pdeg_P f >0$. 
}
\end{definition}

\begin{definition}[cf.\cite{Oka2018},Definition 9.18;\ \cite{Oka2015},p.174]\label{strongly-mixed-wh-face-type}  %% germ
{\rm 
Let $(f,\zero)$ be a mixed function germ at $\zero \in \bc^n$. 
\begin{itemize}
\item The germ $(f,\zero)$ of a mixed function $f(\bbz, \bar{\bbz})$ at $\zero \in \bc^n$ is called 
{\em of mixed weighted homogeneous face type} if %%% あとで、例で使うので追加
the face function $f_\Delta(\bbz, \bar{\bbz})$ is 
a mixed weighted homogeneous polynomial (Definition \ref{mixed-whp}) 
for every \underline{compact} face $\Delta$. 

\item The germ $(f,\zero)$ of a mixed function $f(\bbz, \bar{\bbz})$ at $\zero \in \bc^n$ is called 
{\em of strongly polar positive mixed weighted homogeneous face type} if 
for every 
%%%%%%%%%%%%%%%%%%%%%%%%%%%%%%%%%%%%%%%%%%%%%%%%
\footnote{In \cite{Oka2015}, a convenient (see Definition \ref{convenient-function} of this paper) mixed function $f(\bbz, \bar{\bbz})$ is called %%
{\em of strongly polar positive (mixed) weighted homogeneous face type} if 
the face function $f_\Delta(\bbz, \bar{\bbz})$ is a strongly polar positive (mixed) weighted homogeneous polynomial 
for every $(n-1)$-dimensional face. %%
After this definition, Proposition 10 of \cite{Oka2015} proves that 
for a convenient mixed function $f$ of strongly polar positive (mixed) weighted homogeneous face type (in the sense of \cite{Oka2015}) and 
any weight vector $P$, 
the face function $f_{\Delta(P)}$ is also a strongly polar positive mixed weighted homogeneous polynomial with respect to $P$.}
%%%%%%%%%%%%%%%%%%%%%%%%%%%%%%%%%%%%%%%%%%%%%%%%
compact face $\Delta$, 
the face function $f_\Delta(\bbz, \bar{\bbz})$ is 
a strongly polar positive mixed weighted homogeneous polynomial (Definition \ref{mixed-whp}) 
with respect to 
{\em some} strictly positive weight vector $P$ with $\Delta = \Delta(P)$. %%%
\end{itemize}
}
\end{definition}

\medskip

\section{Newton non-degeneracy and Strong Newton non-degeneracy}\label{section-Newton-non-degeneracy}

\begin{definition}[\cite{Oka2010}, p.6, Definition 3;\ \cite{Oka2018}, p.80 and pp.181--182]\label{Newton-non-degeneracy}
{\rm 
Let $(f,\zero)$ be the germ of a mixed function $f$ at $\zero \in \bc^n$. 
\begin{enumerate}
\item We say $(f,\zero)$ is {\em Newton non-degenerate 
over a compact face $\Delta$} 
if 
$0$ is not a mixed critical value of the face function 
$f_\Delta : {\bc^*}^n \to \bc$. 
(In particular, 
if $f_\Delta^{-1}(0) \cap {\bc^*}^n = \emptyset$, then $0$ is not a mixed critical value of the face function 
$f_\Delta : {\bc^*}^n \to \bc$.) 

\item Let $\Delta$ be a compact face with $\dim \Delta \geq 1$. 
We say $(f,\zero)$ is {\em strongly Newton non-degenerate 
over $\Delta$} 
if 
the face function $f_\Delta : {\bc^*}^n \to \bc$ has no mixed critical points 
and $f_\Delta : {\bc^*}^n \to \bc$ is surjective onto $\bc$. 

\item Let $\Delta$ be a compact face with $\dim \Delta = 0$, that is, 
$\Delta$ is a vertex of $\Gamma_+(f)$. 
We say $(f,\zero)$ is {\em strongly Newton non-degenerate 
over $\Delta$} 
if 
the face function $f_\Delta : {\bc^*}^n \to \bc$ has no mixed critical points. 
(In this case we do not need the surjectivity from ${\bc^*}^n$ to $\bc$.) 
\end{enumerate}
}
\end{definition}

\begin{definition}[\cite{Oka2018}, p.80 and p.182]
{\rm 
We say the germ $(f,\zero)$ of a mixed function $f$ at $\zero \in \bc^n$ is 
{\em Newton non-degenerate} (respectively, {\em strongly Newton non-degenerate}) %%%
if 
$(f,\zero)$ is Newton non-degenerate (respectively, strongly Newton non-degenerate) 
over every compact face $\Delta$. 
}
\end{definition}

\begin{example}[\cite{Oka2018},\ 8.3.1]
The germ of the mixed homogeneous polynomial 
$\rho(\bbz,\overline{\bbz}) = \sum_{j=1}^n z_j \overline{z_j} = \sum_{j=1}^{n}|z_{j}|^{2}$  %%%2次斉次
of degree $2$ at $\zero \in \bc^{n}$
is Newton nondegenerate, but not strongly Newton nondegenerate. 
Actually, for every compact face $\Delta$, we have ${\rho_{\Delta}}^{-1}(0) \cap \bc^{*n} = \emptyset$. 
However, since $\rho_{\Delta}(\bc^{n}) \subset \br$, 
every point in $\bc^{*n}$ is a mixed critical point. Hence, $\rho$ is not strongly Newton nondegenerate. 
\end{example}

Now note that $f : {\bc^*}^n \to \bc$ is surjective 
under the assumption that ``$f^{-1}(0) \cap {\bc^*}^n \neq \emptyset$". 
Namely, in this case, Newton non-degeneracy over a compact face $\Delta(P)$ implies 
strong Newton non-degeneracy over $\Delta(P)$. 
See the next Proposition \ref{remark4} for more precise statements: 
\begin{proposition}[\cite{Oka2010}, Remark 4;\ \cite{Takashimizu2021thesis} and \cite{Saito-Takashimizu2021winter}]\label{remark4}
Let $f(\bbz)$ be a holomorphic weighted homogeneous polynomial 
of positive degree with respect to a strictly positive weight vector $P$. 
(Then, $f = f_{\Delta(P)}$.)

{\bf (i)}\ Suppose that $f$ is Newton non-degenerate over $\Delta(P)$, 
namely, 
$0$ is not a critical value of $f_{\Delta(P)}=f : {\bc^*}^n \to \bc$. 
Then $f_{\Delta(P)}=f : {\bc^*}^n \to \bc$ has no critical point. 
Hence, with {\bf (ii)} below, $f$ is strongly Newton non-degenerate over $\Delta(P)$. 

{\bf (ii)}\ Suppose that $\dim \Delta(P) \geq 1$, namely, $f = f_{\Delta(P)}$ has at least two monomials. 
Then $f_{\Delta(P)}=f : {\bc^*}^n \to \bc$ is surjective. 

\medskip

\noindent 
Let $f(\bbz, \bar{\bbz})$ be a mixed weighted homogeneous polynomial (Definition \ref{mixed-whp}) with respect to 
a radial weight vector $P\ (\gg 0)$ and a polar weight vector $Q$. 

{\bf (iii)}\ Suppose that $(f,\zero)$ is Newton non-degenerate 
\underline{over a compact face $\Delta(P)$} %%% 
and 
the face function $f_{\Delta(P)}=f$ is a polar weighted homogeneous polynomial of \underline{non-zero} polar degree  %%% non-zero polar degree
with respect to the polar weight vector $Q$. 
Then $f : {\bc^*}^n \to \bc$ has no mixed critical points. 

{\bf (iv)}\ In addition to {\bf (iii)}, %%%
we assume that 
$f^{-1}(0) \cap {\bc^*}^n \neq \emptyset$.  %%
Then $f : {\bc^*}^n \to \bc$ is surjective. 
Hence, with {\bf (iii)}, 
in this case, Newton non-degeneracy \underline{over a compact face $\Delta(P)$} implies 
strong Newton non-degeneracy \underline{over $\Delta(P)$}. \hfill \QED
\end{proposition}

The facts in Proposition \ref{remark4} are mentioned in Remark 4 of \cite{Oka2010}. 
We would like to give their detailed proofs in \S \ref{section-appendix}. 

\medskip

\section{The main theorem and true non-degeneracy}\label{section-main}

Related to the assertions {\bf (iii)} and {\bf (iv)} of Proposition \ref{remark4} in the previous section, 
we are interested in the condition 
$$
f^{-1}(0) \cap {\bc^*}^n \neq \emptyset .
$$
Theorem \ref{true-nondegenerate-Thm} below shows 
the existence of a collection of mixed weighted homogeneous polynomials $g(\bbz, \bar{\bbz})$ with respect to 
radial weight vectors $P\ (\gg 0)$ and 
polar weight vectors $Q$ of \underline{non-zero} polar degree 
which satisfy 
$\dim \Delta (P) \geq 1$ and 
$g^{-1}(0) \cap {\bc^*}^n = \emptyset$: 

\begin{theorem}\label{true-nondegenerate-Thm}
Let us consider the mixed polynomial
$$
g(z_{1},z_{2},\overline{z_{1}},\overline{z_{2}}) := 
  z_{1}^{a}\overline{z_{1}}^{\alpha - a}  z_{2}^{b}\overline{z_{2}}^{\beta - b}
- z_{1}^{a'}\overline{z_{1}}^{\alpha - a'}z_{2}^{b'}\overline{z_{2}}^{\beta - b'}
+ z_{1}^{c}\overline{z_{1}}^{\gamma - c}z_{2}^{d}\overline{z_{2}}^{\delta - d}
+ z_{1}^{c'}\overline{z_{1}}^{\gamma - c'}z_{2}^{d'}\overline{z_{2}}^{\delta - d'} ,
$$
where $\alpha, \beta, \gamma$ and $\delta$ are non-negative integers with 
$\alpha < \gamma$ and $\delta < \beta$. 
Let 
$P = {}^t(p_1,p_2)$ be a strictly positive weight vector and 
$Q = {}^t(q_1,q_2) \neq \zero$ be a weight vector $(\in \br_+^2)$. 
We assume that 
\begin{equation}\label{radially-w-h-2}
p_1 (\gamma - \alpha) + p_2 (\delta - \beta) = 0,  %%% 
\end{equation}

\begin{equation}\label{a-a'b-b'}
a-a' = q_2k \ \ \text{and}\ \ b-b' = - q_1k \ \ \ \text{for\ some\ non-zero\ integer}\ k , 
\end{equation}

\begin{equation}\label{condition1}
(a-a',\ b-b') = (c-c',\ d- d') , %%全く等しいとしておく
\end{equation}

\begin{equation}\label{condition2}
(2(a-c)-\alpha + \gamma,\ 2(b-d)-\beta + \delta) = \pm(a-a',\ b-b') , 
\end{equation}
and 
\begin{equation}\label{non-zero-pol-degree}
q_1 (2a - \alpha) + q_2 (2b - \beta) \neq 0 . 
\end{equation}
Then, 
$g$ is a mixed weighted homogeneous polynomial of radial degree $\rdeg_P g = d_r \ (>0)$ with respect to the radial weight vector $P$ and 
of polar degree $\pdeg_Q g = d_p \ (\neq 0)$ with respect to the polar weight vector $Q$. 
We have $\dim \Delta(P)=1$, $g = g_{\Delta(P)}$,     %%% 1次元面に対する面関数
and 
$$
g^{-1}(0)\cap \bc^{*2} = \emptyset . 
$$
Hence, $g$ is Newton non-degenerate and not strongly Newton non-degenerate over $\Delta(P)$. 
\end{theorem}

\begin{proof}
Let us consider the radial Newton polyhedron $\Gamma_{+}(g)$ of the germ $(g,\zero)$. 
Then note that 
the first and second terms correspond to the $0$-dimensional face $\{ (\alpha, \beta) \}$, and 
the third and fourth terms correspond to the $0$-dimensional face $\{ (\gamma, \delta) \}$. %%% １頂点で極重さベクトルが決まってしまう。

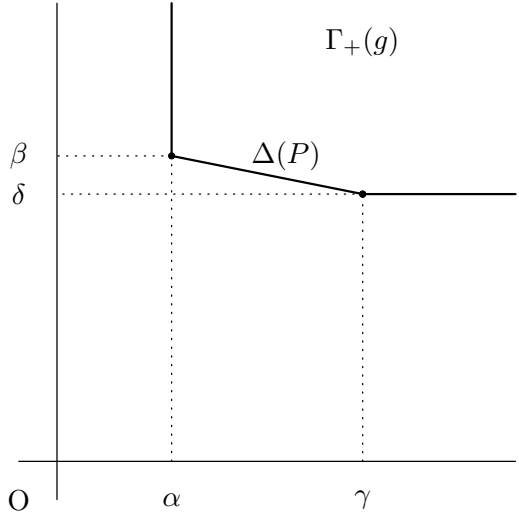
\begin{figure}[H]
\begin{center}
%%%  \input{remark4_example_P_general.tex}
%
%WinTpicVersion4.32a
{\unitlength 0.1in%
\begin{picture}(26.9000,26.0000)(13.1000,-44.0000)%
% LINE 2 2 3 0 Black White  
% 8 2200 4200 2200 2600 2200 2600 1600 2600 3200 4200 3200 2800 3200 2800 1600 2800
% 
\special{pn 8}%
\special{pa 2200 4200}%
\special{pa 2200 2600}%
\special{dt 0.045}%
\special{pa 2200 2600}%
\special{pa 1600 2600}%
\special{dt 0.045}%
\special{pa 3200 4200}%
\special{pa 3200 2800}%
\special{dt 0.045}%
\special{pa 3200 2800}%
\special{pa 1600 2800}%
\special{dt 0.045}%
% DOT 0 0 3 0 Black White  
% 3 2200 2600 3200 2800 3200 2800
% 
\special{pn 4}%
\special{sh 1}%
\special{ar 2200 2600 16 16 0 6.2831853}%
\special{sh 1}%
\special{ar 3200 2800 16 16 0 6.2831853}%
\special{sh 1}%
\special{ar 3200 2800 16 16 0 6.2831853}%
% LINE 1 0 3 0 Black White  
% 6 2200 2600 3200 2800 2200 2600 2200 1800 3200 2800 4000 2800
% 
\special{pn 13}%
\special{pa 2200 2600}%
\special{pa 3200 2800}%
\special{fp}%
\special{pa 2200 2600}%
\special{pa 2200 1800}%
\special{fp}%
\special{pa 3200 2800}%
\special{pa 4000 2800}%
\special{fp}%
% LINE 2 0 3 0 Black White  
% 6 1600 4200 4000 4200 1600 4200 1400 4200 1600 4400 1600 1800
% 
\special{pn 8}%
\special{pa 1600 4200}%
\special{pa 4000 4200}%
\special{fp}%
\special{pa 1600 4200}%
\special{pa 1400 4200}%
\special{fp}%
\special{pa 1600 4400}%
\special{pa 1600 1800}%
\special{fp}%
% STR 2 0 3 0 Black White  
% 4 1400 2500 1400 2600 5 0 0 0
% $8$
\put(14.0000,-26.0000){\makebox(0,0){$\beta$}}%
% STR 2 0 3 0 Black White  
% 4 1400 2700 1400 2800 5 0 0 0
% $7$
\put(14.0000,-28.0000){\makebox(0,0){$\delta$}}%
% STR 2 0 3 0 Black White  
% 4 1400 4300 1400 4400 5 0 0 0
% O
\put(14.0000,-44.0000){\makebox(0,0){O}}%
% STR 2 0 3 0 Black White  
% 4 2200 4300 2200 4400 5 0 0 0
% $3$
\put(22.0000,-44.0000){\makebox(0,0){$\alpha$}}%
% STR 2 0 3 0 Black White  
% 4 3200 4300 3200 4400 5 0 0 0
% $8$
\put(32.0000,-44.0000){\makebox(0,0){$\gamma$}}%
% STR 2 0 3 0 Black White  
% 4 2800 2500 2800 2600 5 0 0 0
% $\Delta(Q)$
\put(28.0000,-26.0000){\makebox(0,0){$\Delta(P)$}}%
% STR 2 0 3 0 Black White  
% 4 3200 1900 3200 2000 5 0 0 0
% $\Gamma_{+}(g)$
\put(32.0000,-20.0000){\makebox(0,0){$\Gamma_{+}(g)$}}%
\end{picture}}%
\end{center}

\vspace{7mm}

\caption{The radial Newton polyhedron of $g(z_{1},z_{2},\overline{z_{1}},\overline{z_{2}})$.}
\label{remark4_example_general}
\end{figure}

By the assumption \eqref{radially-w-h-2} we have 
\begin{equation}\label{radially-w-h}
p_1 \alpha + p_2 \beta = p_1 \gamma + p_2 \delta = d_r \ (>0) .
\end{equation}
Hence, $\dim \Delta(P)=1$ and $g = g_{\Delta(P)}$. 
Let us consider the following condition:
\begin{eqnarray}\label{polar-w-h}
   &q_1 (2a - \alpha) + q_2 (2b - \beta)\\ 
=  &q_1 (2a' - \alpha) + q_2 (2b' - \beta) \nonumber \\
=  &q_1 (2c - \gamma) + q_2 (2d - \delta) \ \nonumber \\
=  &q_1 (2c' - \gamma) + q_2 (2d' - \delta) \nonumber \\
=  &d_p \ (\neq 0).\ \ \ \ \ \ \ \ \ \ \ \ \ \ \ \ \ \ \ \nonumber 
\end{eqnarray}
By \eqref{a-a'b-b'} we have 
\begin{equation}\label{aba'b'}
  q_1 (2a - \alpha) + q_2 (2b - \beta) 
= q_1 (2a' - \alpha) + q_2 (2b' - \beta) .
\end{equation}
By \eqref{condition1} we have 
\begin{equation}\label{cdc'd'}
  q_1 (2c - \gamma) + q_2 (2d - \delta) 
= q_1 (2c' - \gamma) + q_2 (2d' - \delta). 
\end{equation}
By \eqref{condition2}  we have 
$$
q_1(2a-2c-\alpha + \gamma) + q_2(2b-2d-\beta + \delta) = 0 ,
$$
and hence, 
\begin{equation}\label{abcd}
  q_1(2a - \alpha) + q_2(2b - \beta) 
= q_1(2c - \gamma) + q_2(2d - \delta) .
\end{equation}
%%%%%%%%
We see that \eqref{aba'b'}, \eqref{cdc'd'}, \eqref{abcd} and \eqref{non-zero-pol-degree}
imply \eqref{polar-w-h}. 
%%%%%%%%
From the conditions \eqref{radially-w-h} and \eqref{polar-w-h}, 
$g$ is a mixed weighted homogeneous polynomial 
of radial degree $\rdeg_P g = d_r \ (>0)$ 
with respect to the radial weight vector $P$ and 
of polar degree $\pdeg_Q g = d_p \ (\neq 0)$ 
with respect to the polar weight vector $Q$. 
From \eqref{condition1} and \eqref{condition2} we have 
\begin{equation}\label{difference}
((a+a'-\alpha)-(c+c'-\gamma), (b+b'-\beta)-(d+d'-\delta)) 
%%%% = (a-c + a'- c'-\alpha + \gamma, b-d + b'- d'-\beta + \delta) 
=  \pm(a-a', b-b') .
\end{equation}
We see that $g = g_{\Delta(P)}$ and $\dim \Delta(P)=1$. %%% 1次元面に対する面関数
We will show that 
$g^{-1}(0)\cap \bc^{*2} = \emptyset$. 

For $(z_{1},z_{2}) \in \bc^{*2}$ with $z_{1}=r_{1}e^{i\theta_{1}},z_{2}=r_{2}e^{i\theta_{2}}\ \ (r_{1}>0,\ r_{2}>0)$, we have 
$$
\begin{array}{l}
g(z_{1},z_{2},\overline{z_{1}},\overline{z_{2}})\\
=
  r_{1}^{\alpha}r_{2}^{\beta}e^{i((2a - \alpha)\theta_{1}+ (2b - \beta)\theta_{2})}
- r_{1}^{\alpha}r_{2}^{\beta}e^{i((2a' - \alpha)\theta_{1}+ (2b' - \beta)\theta_{2})}
+ r_{1}^{\gamma}r_{2}^{\delta}e^{i((2c - \gamma)\theta_{1}+ (2d - \gamma)\theta_{2})}
+ r_{1}^{\gamma}r_{2}^{\delta}e^{i((2c' - \gamma)\theta_{1}+ (2d' - \gamma)\theta_{2})}\\
=
 r_{1}^{\alpha}r_{2}^{\delta}
(r_{2}^{\beta - \delta}e^{i((2a - \alpha)\theta_{1}+ (2b - \beta)\theta_{2})}
 - r_{2}^{\beta - \delta}e^{i((2a' - \alpha)\theta_{1}+ (2b' - \beta)\theta_{2})}
 + r_{1}^{\gamma - \alpha}e^{i((2c - \gamma)\theta_{1}+ (2d - \gamma)\theta_{2})}\\
\hspace{6cm} + r_{1}^{\gamma - \alpha}e^{i((2c' - \gamma)\theta_{1}+ (2d' - \gamma)\theta_{2})}) \\
=
 r_{1}^{\alpha}r_{2}^{\delta}
( r_{2}^{\beta - \delta}(e^{i((2a - \alpha)\theta_{1}+ (2b - \beta)\theta_{2})}
 - e^{i((2a' - \alpha)\theta_{1}+ (2b' - \beta)\theta_{2})})\\
\hspace{6cm} + r_{1}^{\gamma - \alpha}(e^{i((2c - \gamma)\theta_{1}+ (2d - \gamma)\theta_{2})}
 + e^{i((2c' - \gamma)\theta_{1}+ (2d' - \gamma)\theta_{2})}) ) .
\end{array}
$$
Here we have 
$$
\begin{array}{cl}
& \left( r_{2}^{\beta - \delta}(e^{i((2a - \alpha)\theta_{1}+ (2b - \beta)\theta_{2})}
   - e^{i((2a' - \alpha)\theta_{1}+ (2b' - \beta)\theta_{2})})
   + r_{1}^{\gamma - \alpha}(e^{i((2c - \gamma)\theta_{1}+ (2d - \gamma)\theta_{2})}
   + e^{i((2c' - \gamma)\theta_{1}+ (2d' - \gamma)\theta_{2})}) \right)\\
= &
2\,i\,r_{2}^{\beta - \delta}\sin \frac{1}{2}\left( ((2a - \alpha)\theta_{1}+ (2b - \beta)\theta_{2}) - ((2a' - \alpha)\theta_{1}+ (2b' - \beta)\theta_{2}) \right)\\
   &\times e^{i (((2a - \alpha)\theta_{1}+ (2b - \beta)\theta_{2})+((2a' - \alpha)\theta_{1}+ (2b' - \beta)\theta_{2}))/2}\\
& +\, 2\,r_{1}^{\gamma - \alpha}\cos \frac{1}{2}\left( ((2c - \gamma)\theta_{1}+ (2d - \delta)\theta_{2}) - ((2c' - \gamma)\theta_{1}+ (2d' - \delta)\theta_{2}) \right)\\
   &\times e^{i (((2c - \gamma)\theta_{1}+ (2d - \delta)\theta_{2})+((2c' - \gamma)\theta_{1}+ (2d' - \delta)\theta_{2}))/2}\\
= &
2\,i\,r_{2}^{\beta - \delta}\sin ((a -a')\theta_{1}+ (b - b')\theta_{2}) e^{i((a + a' - \alpha)\theta_{1}+ (b + b' - \beta)\theta_{2})}\\
& +\, 2\,r_{1}^{\gamma - \alpha}\cos ((c - c')\theta_{1}+ (d - d')\theta_{2}) e^{i((c + c' - \gamma)\theta_{1}+ (d + d' - \delta)\theta_{2})}\\
= &
2( r_{2}^{\beta - \delta}\sin ((a -a')\theta_{1}+ (b - b')\theta_{2}) e^{i((a + a' - \alpha)\theta_{1}+ (b + b' - \beta)\theta_{2} + \pi/2)} \\
& +\, r_{1}^{\gamma - \alpha}\cos ((c - c')\theta_{1}+ (d - d')\theta_{2}) e^{i((c + c' - \gamma)\theta_{1}+ (d + d' - \delta)\theta_{2})} ) \\
= &
2( r_{2}^{\beta - \delta}\sin ((a -a')\theta_{1}+ (b - b')\theta_{2}) e^{i((a + a' - \alpha)\theta_{1}+ (b + b' - \beta)\theta_{2} + \pi/2)} \\
& +\, r_{1}^{\gamma - \alpha}\cos ((a - a')\theta_{1}+ (b - b')\theta_{2}) e^{i((c + c' - \gamma)\theta_{1}+ (d + d' - \delta)\theta_{2})} )  .
\end{array}
$$
If one of $\sin ((a -a')\theta_{1}+ (b - b')\theta_{2}),\ \cos ((a - a')\theta_{1}+ (b - b')\theta_{2})$ is $0$, then the other is not $0$, 
and hence, we have $g(z_{1},z_{2},\overline{z_{1}},\overline{z_{2}})\neq 0$. 
Thus, 
if $g(z_{1},z_{2},\overline{z_{1}},\overline{z_{2}})=0$, then we have 
$\sin ((a -a')\theta_{1}+ (b - b')\theta_{2}) \neq 0$ and $\cos ((a - a')\theta_{1}+ (b - b')\theta_{2}) \neq 0$. 
From the equality 
$r_{2}^{\beta - \delta}\sin ((a -a')\theta_{1}+ (b - b')\theta_{2}) e^{i((a + a' - \alpha)\theta_{1}+ (b + b' - \beta)\theta_{2} + \pi/2)} 
 + 
r_{1}^{\gamma - \alpha}\cos ((a - a')\theta_{1}+ (b - b')\theta_{2}) e^{i((c + c' - \gamma)\theta_{1}+ (d + d' - \delta)\theta_{2})}
  =  0$, 
%%%% 両辺を\cos で割ると
we have 
$$
r_{2}^{\beta - \delta}\tan ((a -a')\theta_{1}+ (b - b')\theta_{2}) e^{i((a + a' - \alpha)\theta_{1}+ (b + b' - \beta)\theta_{2} + \pi/2)} 
 + r_{1}^{\gamma - \alpha} e^{i((c + c' - \gamma)\theta_{1}+ (d + d' - \delta)\theta_{2})}  =  0
.$$
Thus, we have 
$$
r_{2}^{\beta - \delta}\tan ((a -a')\theta_{1}+ (b - b')\theta_{2}) 
e^{i( ((a + a' - \alpha) - (c + c' - \gamma) )\theta_{1} +  ((b + b' - \beta) - (d + d' - \delta) )\theta_{2} + \pi/2)} 
 + r_{1}^{\gamma - \alpha}  =  0 .
$$
Using \eqref{difference}, we have 
$$
r_{1}^{\gamma - \alpha} = 
- r_{2}^{\beta - \delta}\tan ((a -a')\theta_{1}+ (b - b')\theta_{2})e^{i(\pm ((a-a')\theta_{1} +  (b-b')\theta_{2}) + \pi/2)} .
$$
Hence, $e^{i(\pm ((a-a')\theta_{1} +  (b-b')\theta_{2}) + \pi/2)}$ is a real number, and 
$$
\pm ((a-a')\theta_{1} +  (b-b')\theta_{2}) + \pi/2 = m\pi
$$
for some integer $m$. 
Namely, 
$$
((a-a')\theta_{1} +  (b-b')\theta_{2}) \pm \pi/2 = m'\pi
$$
for some integer $m'$. 
Then we have $\cos ((a-a')\theta_{1} +  (b-b')\theta_{2}) = 0$. 
This is a contradiction. 
Thus we have $g^{-1}(0)\cap \bc^{*2} = \emptyset$. 
\end{proof}
%%%%%%%%%%%%%%%%%%%%%%%%%%%%%%%%%%%%%%%%%%%%%%%%%%

\begin{example}\label{example3887}
As an example of Theorem \ref{true-nondegenerate-Thm}, 
let us consider the mixed polynomial
$$
g(z_{1},z_{2},\overline{z_{1}},\overline{z_{2}}) := 
  z_{1}^{a}\overline{z_{1}}^{3 - a}  z_{2}^{b}\overline{z_{2}}^{8 - b}
- z_{1}^{a'}\overline{z_{1}}^{3 - a'}z_{2}^{b'}\overline{z_{2}}^{8 - b'}
+ z_{1}^{c}\overline{z_{1}}^{8 - c}z_{2}^{d}\overline{z_{2}}^{7 - d}
+ z_{1}^{c'}\overline{z_{1}}^{8 - c'}z_{2}^{d'}\overline{z_{2}}^{7 - d'} .
$$
%%%%
Then $g$ is a radially weighted homogeneous polynomial of radial degree $43$ 
with respect to the radial weight vector $P := {}^t(1, 5)$. 
We set $Q := {}^t(1, 1)$. 
Using Theorem \ref{true-nondegenerate-Thm}, 
let us enumerate up all 
$$(a,b,a',b',c,d,c',d')\text{'s}$$
satisfying the conditions: 
\begin{equation}\label{condition_I}
(a-a',\ b-b') = (c-c',\ d- d') = (-1,\ 1) , %%
\end{equation}
\begin{equation}\label{condition_II}
(2(a-c)-3 + 8,\ 2(b-d)-8 + 7) = \pm(-1,\ 1) , 
\end{equation}
and 
\begin{equation}\label{positive-pol-degree}
(2a - 3) + (2b - 8) > 0 . %%%%正のものだけ考える
\end{equation}
Here we restrict ourselves to the cases where $g$ is a polar weighted homogeneous polynomial 
of positive polar degree $2a + 2b - 11$ with respect to the polar weight vector $Q = {}^t(1, 1)$ 
by the condition \eqref{positive-pol-degree}. 

The condition \eqref{positive-pol-degree} is equivalent to 
$
a+b \geq 6 .
$
By \eqref{condition_I}, we have $a' = a+1,\ b' = b-1,\ c' = c+1,\ \text{and}\ d' = d-1$. 
Hence, $(a', b', c', d')$ is determined by $(a, b, c, d)$, and 
$
0\leq a \leq 2,\ 1\leq b \leq 8,\ 0\leq c \leq 7,\ \text{and}\ 1\leq d \leq 7 .
$
We have 
$$(2(a-c)-3 + 8,\ 2(b-d)-8 + 7) = (-1,\ 1) \Longleftrightarrow \text{I:}\ \ \ c=a+3\ \text{and}\ d=b-1 ,$$
and 
$$(2(a-c)-3 + 8,\ 2(b-d)-8 + 7) = (1,\ -1) \Longleftrightarrow \text{II:}\ \ \ c=a+2\ \text{and}\ d=b .$$
Finally, $(c, d)$ is determined by $(a, b)$. 
In Case I, $(a,b) \in \bz^2$ should satisfy the conditions
$$
a+b \geq 6,\ 0\leq a \leq 2,\ 2\leq b \leq 8 ,
$$
and 
in Case II, $(a,b) \in \bz^2$ should satisfy the conditions
$$
a+b \geq 6,\ 0\leq a \leq 2,\ 1\leq b \leq 7 .
$$
See Figure \ref{ab}. 
We see that there are 21 $(a,b,a',b',c,d,c',d')$\, 's (mixed weighted homogeneous polynomials). 
\end{example}

\begin{figure}[H]
\begin{center}
\begin{picture}(90,100)
\put(0,0){\vector(0,1){82}}
\put(0,0){\vector(1,0){72}}
\put( -2,86){{\footnotesize $b$}} %vertical
\put(76, -2){{\footnotesize $a$}} %horizontal

\multiput(8,0)(8,0){7}{\line(0,1){80}}
\multiput(0,8)(0,8){9}{\line(1,0){70}}

\put(16,-22){\tiny Case I.}

\put( -2,-10){{\tiny $0$}}
\put(  6,-10){{\tiny $1$}}
\put( 14,-10){{\tiny $2$}}
\put( 22,-10){{\tiny $3$}}
\put( 30,-10){{\tiny $4$}}
\put( 38,-10){{\tiny $5$}}
\put( 46,-10){{\tiny $6$}}
\put( 54,-10){{\tiny $7$}}

\put(-8, -1){{\tiny $0$}}
\put(-8,  7){{\tiny $1$}}
\put(-8, 15){{\tiny $2$}}
\put(-8, 23){{\tiny $3$}}
\put(-8, 31){{\tiny $4$}}
\put(-8, 39){{\tiny $5$}}
\put(-8, 47){{\tiny $6$}}
\put(-8, 55){{\tiny $7$}}
\put(-8, 63){{\tiny $8$}}
\put(-8, 71){{\tiny $9$}}

\put(16,32){\circle*{3}}
\put( 8,40){\circle*{3}}
\put(16,40){\circle*{3}}
\put( 0,48){\circle*{3}}
\put( 8,48){\circle*{3}}
\put(16,48){\circle*{3}}
\put( 0,56){\circle*{3}}
\put( 8,56){\circle*{3}}
\put(16,56){\circle*{3}}
\put( 0,64){\circle*{3}}
\put( 8,64){\circle*{3}}
\put(16,64){\circle*{3}}
\end{picture}
\ \ \ \ \ \ \ \ \ 
\begin{picture}(90,100)
\put(0,0){\vector(0,1){82}}
\put(0,0){\vector(1,0){72}}
\put( -2,86){{\footnotesize $b$}} %vertical
\put(76, -2){{\footnotesize $a$}} %horizontal

\multiput(8,0)(8,0){7}{\line(0,1){80}}
\multiput(0,8)(0,8){9}{\line(1,0){70}}

\put(16,-22){\tiny Case II.}

\put( -2,-10){{\tiny $0$}}
\put(  6,-10){{\tiny $1$}}
\put( 14,-10){{\tiny $2$}}
\put( 22,-10){{\tiny $3$}}
\put( 30,-10){{\tiny $4$}}
\put( 38,-10){{\tiny $5$}}
\put( 46,-10){{\tiny $6$}}
\put( 54,-10){{\tiny $7$}}

\put(-8, -1){{\tiny $0$}}
\put(-8,  7){{\tiny $1$}}
\put(-8, 15){{\tiny $2$}}
\put(-8, 23){{\tiny $3$}}
\put(-8, 31){{\tiny $4$}}
\put(-8, 39){{\tiny $5$}}
\put(-8, 47){{\tiny $6$}}
\put(-8, 55){{\tiny $7$}}
\put(-8, 63){{\tiny $8$}}
\put(-8, 71){{\tiny $9$}}

\put(16,32){\circle*{3}}
\put( 8,40){\circle*{3}}
\put(16,40){\circle*{3}}
\put( 0,48){\circle*{3}}
\put( 8,48){\circle*{3}}
\put(16,48){\circle*{3}}
\put( 0,56){\circle*{3}}
\put( 8,56){\circle*{3}}
\put(16,56){\circle*{3}}
\end{picture}
\end{center}
\vspace*{10mm}
\caption{\ \ }
\label{ab}
\end{figure}

\begin{example}\label{example3967}
Let us consider the mixed polynomial
$$
h(z_{1},z_{2},\overline{z_{1}},\overline{z_{2}}) := 
  \overline{z_{1}}^{3}  z_{2}^{6}\overline{z_{2}}^{3}
- z_{1}^{3}z_{2}^{4}\overline{z_{2}}^{5}
+ \overline{z_{1}}^{6}z_{2}^{6}\overline{z_{2}}
+ z_{1}^{3}\overline{z_{1}}^{3}z_{2}^{4}\overline{z_{2}}^{3} .
$$
%%%
Then $h$ is a \underline{strongly} mixed weighted homogeneous polynomial of 
radial degree $33$ and polar degree $3$ 
with respect to the weight vector $P := {}^t(2, 3)$. 
We have $\dim \Delta(P)=1$ and $h = h_{\Delta(P)}$.      %%% 1次元面に対する面関数
Since $h$ satisfies the conditions 
$
(a-a',\ b-b') = (c-c',\ d- d') = (-3,\ 2) %%全く等しい
$
and 
$
(2(a-c)-\alpha + \gamma,\ 2(b-d)-\beta + \delta) = (3,\ -2)
$, 
by Theorem \ref{true-nondegenerate-Thm}, 
we have 
$
h^{-1}(0)\cap \bc^{*2} = \emptyset . 
$
\end{example}

\bigskip

Related to the condition $f^{-1}(0) \cap {\bc^*}^n \neq \emptyset$ in Proposition \ref{remark4}, 
the notion ``true non-degeneracy" is introduced in Oka \cite{Oka2010}: 
\begin{definition}[\cite{Oka2010}, Definition 3]\label{true-non-degenerate}
{\rm 
A mixed function germ $f(\bbz,\overline{\bbz})$ at $\zero \in \bc^n$ is called 
{\em true non-degenerate} if 
$f$ is Newton non-degenerate and 
$f_{\Delta}^{-1}(0) \cap {\bc^*}^n \neq \emptyset$ 
for every compact face $\Delta$ of $\Gamma_+(f)$ with $\dim \Delta \geq 1$. 
}
\end{definition}

Now let $(f,\zero)$ be the germ of a mixed funtion 
$f(\bbz,\overline{\bbz})=\sum_{\nu,\mu}c_{\nu,\mu}\bbz^{\nu}\overline{\bbz}^{\mu}$. 
We assume that $f(\zero)=0$. 
For a subset $I \subset \{1,2,\cdots,n\}$, we set 
$\bc^I := \{ \bbz \,|\, z_j= 0\ \text{for\ every}\ j\notin I \}$ and 
$f^I :=f|_{\bc^I}$.           %%% fの定義域Uは省略

\begin{definition}\label{convenient-function} %%%関数の芽 f が利便であるとは
{\rm 
A mixed function germ $(f,\zero)$ is {\em convenient} 
if 
$f^{I} \not\equiv 0$ for every $I \subset \{ 1,2,\cdots,n \}$ with $I \neq \emptyset$. 
}
\end{definition}

Note that 
if $(f,\zero)$ is a convenient Newton non-degenerate mixed function germ, then 
%%% Oka（\cite{Oka2010}）のTheorem 19において，主張(1)のみを示すには，true non-degenerateという仮定は不要で，Newton非退化でさえあればよい．
$\zero$ is a mixed regular point of $f$ or an isolated mixed critical point of $f$         %%% convenientならば、高々孤立混合臨界点である．
in the mixed hypersurface $V := f^{-1}(0)$ 
(\cite{Oka2010}, Theorem 19\ (1)). 

\begin{example}
Let us consider the mixed function germ 
$$
f(z_{1},z_{2},\overline{z_{1}},\overline{z_{2}}) := 
z_{2}^{7}\overline{z_{2}}^{2}
+  z_{1}^{3} z_{2}^{6}\overline{z_{2}}^{2}
- \overline{z_{1}}^{3}z_{2}^{7}\overline{z_{2}}
+ z_{1}^{7}\overline{z_{1}}z_{2}^{5}\overline{z_{2}}^{2}
+ z_{1}^{4}\overline{z_{1}}^{4}z_{2}^{6}\overline{z_{2}}
+ z_{1}^{35}\overline{z_{1}}^{20} .
$$
It is obvious that $(f,\zero)$ is convenient
\footnote{
A mixed function germ $(f,\zero)$ is convenient if and only if 
for every fixed $i =1, \dots , n$, there exists a pair of non-negative integers $(\nu_i,\ \mu_i)$ with $\nu_i + \mu_i >0$ 
such that $f$ has the term $c_{\nu, \mu}z_i^{\nu_i}\bar{z}_i^{\mu_i} \ (c_{\nu, \mu}\neq 0)$, 
where $\nu := (0, \dots , 0, \nu_i , 0, \dots , 0),\ \mu := (0, \dots , 0, \mu_i , 0, \dots , 0)$. 
See Lemma \ref{convenient-function-criterion} in \S \ref{section-appendix}.}
. 
We set 
$
P:= {}^t(1, 3),\    Q:= {}^t(1, 5),\ R:= {}^t(7, 47)
$. 
Then we have: 
\begin{itemize}
\item 
The $1$-dimensional face function 
$
f_{\Delta(P)} = z_{2}^{7}\overline{z_{2}}^{2}
+  z_{1}^{3} z_{2}^{6}\overline{z_{2}}^{2}
- \overline{z_{1}}^{3}z_{2}^{7}\overline{z_{2}}
$
is a strongly mixed weighted homogeneous polynomial of radial degree $27$ and polar degree $15$  %%
with respect to the weight vector $P$. 

\item 
The $1$-dimensional face function 
$
f_{\Delta(Q)} = z_{1}^{3} z_{2}^{6}\overline{z_{2}}^{2}
- \overline{z_{1}}^{3}z_{2}^{7}\overline{z_{2}}
+ z_{1}^{7}\overline{z_{1}}z_{2}^{5}\overline{z_{2}}^{2}
+ z_{1}^{4}\overline{z_{1}}^{4}z_{2}^{6}\overline{z_{2}}
$
is a mixed weighted homogeneous polynomial   %%
of radial degree $43$ 
with respect to the radial weight vector $Q$ and 
of polar degree $15$ 
with respect to the polar weight vector $P$. 
By Therem \ref{true-nondegenerate-Thm}, 
we have 
$$
f_{\Delta(Q)}^{-1}(0)\cap \bc^{*2} = \emptyset . 
$$

\item 
The $1$-dimensional face function 
$
f_{\Delta(R)} = z_{1}^{7}\overline{z_{1}}z_{2}^{5}\overline{z_{2}}^{2}
+ z_{1}^{4}\overline{z_{1}}^{4}z_{2}^{6}\overline{z_{2}}
+ z_{1}^{35}\overline{z_{1}}^{20}
$
is a mixed weighted homogeneous polynomial   %%
of radial degree $385$ 
with respect to the radial weight vector $R$ and 
of polar degree $15$ 
with respect to the polar weight vector $P$. 
\end{itemize}
Thus $(f,\zero)$ is a convenient mixed function germ 
\underline{of mixed weighted homogeneous face type} (Definition \ref{strongly-mixed-wh-face-type}) 
which is not true non-degenerate. 
\end{example}

\section{Appendix}\label{section-appendix}

In this section we give some proofs of facts stated in \S \ref{section-mixed-function}--\S \ref{section-main} for unfamiliar readers. 

\begin{lemma}[cf. \cite{Saito-Takashimizu2021-8}]\label{compact-face}
A face $\Delta$ of $\Gamma_+(f)$ is compact 
if and only if 
$\Delta = \Delta(P)$ 
for some strictly positive weight vector $P$. 
\end{lemma}
\begin{proof}
If $P={}^{t}\!(p_{1},\dots,p_{n})$ is strictly positive, then 
the set 
$
\{ (x_1, \dots , x_n) \in \br_+^n \ |\ P(x_1, \dots , x_n) = d(P) \}
$
is an $(n-1)$-simplex whose vertices are 
$
d(P)(1/p_1, 0 , \dots , 0), \dots , d(P)(0 , \dots , 0, 1/p_n)
$. 
Since the face $\Delta(P)$ is contained in this $(n-1)$-simplex, 
it is a bounded closed subset of $\br_+^n$, and hence, it is compact. 
Conversely, if a face $\Delta = \Delta(Q)$ is compact, where $Q = {}^{t}\!(q_1,\dots,q_n)$, 
then we can take $Q$ to be strictly positive. 
Actually, take the minimum number $m \ (1 \leq m \leq n)$ 
such that $\Delta \subset \br^m \times \{ \zero \}$ 
changing the order of the components of $\br^n$ if necessary. 
Then, we have $Q(\bx) = {}^{t}\!(q_1,\dots,q_m , M , \dots , M)(\bx)$ 
for every $\bx \in \Delta$ and every positive real number $M$. 
If we take a sufficiently large positive integer $M_0$, %%% 整数
then we have 
$
Q(\bx) \leq {}^{t}\!(q_1,\dots,q_m , M_0 , \dots , M_0)(\bx)
$
for all $\bx \in \Gamma_{+}(f)$. 
Put $P := {}^{t}\!(q_1,\dots,q_m , M_0 , \dots , M_0)$. Then we have 
$
\Delta = \Delta(Q) = \Delta(P) .
$
By the definition of $m$, 
there exists a point $\bx \in \Delta$ such that $x_i \neq 0$ for every $i \ (1 \leq i \leq m)$. 
Suppose that $q_i = 0$. 
Then $\bx' := (x_1, \dots , x_{i-1}, \alpha , x_{i+1}, \dots ,x_n)\ (\in \Gamma_{+}(f))$,   %%%% $x_i \neq 0$の非負実数倍
where $\alpha$ is an arbitrary non-negative real number, 
is contained in $\Delta(P) = \Delta$. 
This contradicts the compactness of $\Delta$. 
Hence, $P$ is strictly positive. 
\end{proof}

\medskip

We next give the proof of Lemma \ref{mixed-action} as follows: 
If $f(t \circ \bbz) = t^{d_r}f(\bbz)\ \ \ (\forall t>0, \forall \bbz \in \bc^n)$, then 
we have 
$
\sum_{\nu, \mu}t^{P(\nu + \mu)}c_{\nu, \mu} \bbz^\nu \bar{\bbz}^\mu
=
\sum_{\nu, \mu}t^{d_r}c_{\nu, \mu} \bbz^\nu \bar{\bbz}^\mu
$. 
By the uniqueness of the Talor expansion, we have 
$$
c_{\nu, \mu} \neq 0 \Longrightarrow t^{P(\nu + \mu)} = t^{d_r}\ \ \ (\forall t>0). 
$$
Hence, we have $P(\nu + \mu) = d_r$. 
The proof of the converse is easy. 
Now if $f(\theta \circ \bbz) = e^{id_p \theta}f(\bbz)\ \ \ (\forall \theta \in \br, \forall \bbz \in \bc^n)$, then 
we have 
$
\sum_{\nu, \mu}e^{iQ(\nu - \mu)\theta}c_{\nu, \mu} \bbz^\nu \bar{\bbz}^\mu
=
\sum_{\nu, \mu}e^{id_p\theta}c_{\nu, \mu} \bbz^\nu \bar{\bbz}^\mu
$. 
By the uniqueness of the Talor expansion, we have 
$$
c_{\nu, \mu} \neq 0 \Longrightarrow e^{iQ(\nu - \mu)\theta} = e^{id_p\theta}\ \ \ (\forall \theta \in \br).
$$
Take a real number $\theta \ (\neq 0)$ such that 
$
-\pi < Q(\nu - \mu)\theta ,\ d_p\theta \leq \pi .
$
Then we have $Q(\nu - \mu) = d_p$. 
The proof of the converse is easy. \hfill \QED
%%%%%%%%%%%%%%%%%%%%%%%%%%%%%%%%%%%%%%%%%%%%%%%%%%%%%%%%%%%%

\bigskip

We now give a direct proof of Proposition \ref{Euler}. This proof is different from p.183 of \cite{Oka2018}. 
Let us prove {\bf (R)}. If $f$ is a radially weighted homogeneous polynomial 
of radial degree $d_{r}\ (>0)$ with respect to a weight vector $P={}^{t}\!(p_{1},\dots,p_{n})$, then we have 
$$
\frac{\partial f}{\partial z_j}
= \sum_{\nu,\mu}c_{\nu,\mu}\nu_j z_1^{\nu_1} \cdots z_j^{\nu_j -1} \cdots z_n^{\nu_n}\overline{\bbz}^{\mu} .
$$
Hence, we have 
$
p_j z_j \frac{\partial f}{\partial z_j}
= \sum_{\nu,\mu}c_{\nu,\mu}p_j \nu_j z_1^{\nu_1} \cdots z_n^{\nu_n}\overline{\bbz}^{\mu}
$. 
Hence, we have 
$$
\sum_j^n p_j z_j \frac{\partial f}{\partial z_j}
= \sum_{\nu,\mu}c_{\nu,\mu} \left( \sum_j^n p_j \nu_j \bbz^{\nu}\overline{\bbz}^{\mu} \right)
$$
and 
$$
\sum_j^n p_j \overline{z_j} \frac{\partial f}{\partial \overline{z_j}}
= \sum_{\nu,\mu}c_{\nu,\mu} \left( \sum_j^n p_j \mu_j \bbz^{\nu}\overline{\bbz}^{\mu} \right) .
$$
Hence, we have 
$$
\sum_{j=1}^{n}p_{j}\left(z_{j}\frac{\partial f}{\partial z_{j}} + \overline{z_{j}}\frac{\partial f}{\partial \overline{z_{j}}} \right)
= 
\sum_{\nu,\mu}c_{\nu,\mu} \left( \sum_j^n p_j (\nu_j + \mu_j)\right) \bbz^{\nu}\overline{\bbz}^{\mu} = d_{r}f(\bbz,\overline{\bbz}) .
$$
We can prove {\bf (P)} in the same way.  \hfill \QED  %%%% ここで，極次数$d_{p}$は$0$や負の数でもよい．

%%%%%%%%%%%%%%%%%%%%%%%%%%%%%%%%%%%%%%%%%%%%%%%%%%%%%%%%%%%%%%%%%%%%%%%%%%%%%%%%

\bigskip

Here we would like to give a detailed proof of Proposition \ref{remark4}. 
The facts in Proposition \ref{remark4} are mentioned in Remark 4 of \cite{Oka2010}. 

{\bf (i)}\ Suppose that $f$ is Newton non-degenerate over $\Delta(P)$, 
namely, 
$0$ is not a critical value of $f_{\Delta(P)}=f : {\bc^*}^n \to \bc$. 
By Euler equality \eqref{reulereq}, we have 
$$
\sum_{j=1}^{n}p_{j}z_{j}\frac{\partial f}{\partial z_{j}}
= 
(\deg_{P} f) \cdot f(\bbz) .
$$
Hence, if $\bbz \in {\bc^*}^n$ is a critical point of $f$, then we have 
($\deg_{P} f) \cdot f(\bbz) = 0$. 
Since $\deg_{P} f>0$, we have $f(\bbz)=0$. This is a contradiction. 

\medskip

{\bf (ii)}\ Suppose that $\dim \Delta(P) \geq 1$, namely, $f = f_{\Delta(P)}$ has at least two monomials. 
We will show that 
$f_{\Delta(P)}=f : {\bc^*}^n \to \bc$ is surjective by the induction on $n$. 
Here we have $n\geq 2$ since $f = f_{\Delta(P)}$ has at least two monomials. 

We put $f(\bbz) = \sum_{\nu \in \Delta (P)} c_{\nu}\bbz^{\nu},\ 
\bbz=(z_{1},z_{2},\dots,z_{n}),\ 
\nu=(\nu_{1},\nu_{2},\dots,\nu_{n})$, where 
$P = {}^t(p_1,p_2,\dots,p_n) \gg 0$ and $d := \deg_{P} f \ (> 0)$. 

\noindent
Step 1:\ \ Let us consider the case $n=2$. Let $m_{1} >0$ be the highest exponent of $z_{1}$. We have 
$$
f(\bbz)=c_{m_{1},m_{2}}z_{2}^{m_{2}}z_{1}^{m_{1}}+\cdots+c_{n_{1},n_{2}}z_{2}^{n_{2}}z_{1}^{n_{1}},
\ c_{m_{1},m_{2}} \neq 0,\ c_{n_{1},n_{2}} \neq 0 ,
$$
where $n_{1}$ is the lowest exponent of $z_{1}$. 
Since $f$ has at least two monomials, we have $m_{1} > n_{1} \geq 0$. Then we have 
\begin{eqnarray*}
f(z_{1},1) &=& c_{m_{1},m_{2}}z_{1}^{m_{1}}+\cdots+c_{n_{1},n_{2}}z_{1}^{n_{1}} \\
&=& z_{1}^{n_{1}}(c_{m_{1},m_{2}}z_{1}^{m_{1}-n_{1}}+\cdots+c_{n_{1},n_{2}}) .
\end{eqnarray*}
Suppose that $f(z_{1},z_{2}) \equiv 0$ on $\bc^{*2}$. 
Then $f(z_{1},1) = 0$ for every $z_{1} \in \bc^{*}$ and hence, 
there exist infinitely many roots of the equation $f(z_{1},1)=0$. 
Since $f(z_{1},1)$ is a holomorphic polynomial with complex coefficients and $1$ variable of degree $m_1 \ (>0)$, 
this is a contradiction. 
Thus, there exists $(\alpha_{1},\alpha_{2}) \in \bc^{*2}$ such that $f(\alpha_{1},\alpha_{2}) \neq 0$. 
Since $f$ is a weighted homogeneous polynomial of degree $d\ (>0)$ with respect to the weight vector$P$ 
we have 
$$
f(t^{p_{1}}\alpha_{1},t^{p_{2}}\alpha_{2}) = t^{d}f(\alpha_{1},\alpha_{2})
$$
for every $t \in \bc^{*}$. 
Hence, we see that $f(\bc^{*2}) \supset \bc^{*}$. 
%%%% $d>0$より，$t^{d}f(\alpha_{1},\alpha_{2})$は$\bc^{*}$の任意の点を取り得る．
Moreover, we have 
$f(z_{1},1)=0$ if and only if 
$z_{1}=0$ or $c_{m_{1},m_{2}}z_{1}^{m_{1}-n_{1}}+\cdots+c_{n_{1},n_{2}}=0$, where 
$m_{1}-n_{1}>0,\ c_{m_{1},m_{2}} \neq 0,\ c_{n_{1},n_{2}}\neq 0$. 
Hence, the equation $f(z_{1},1)=0$ has a root $z_{1}\ (\neq 0)$. 
Thus, $f:\bc^{*2} \to \bc$ is surjective. 

\noindent
Step 2:\ \ 
Let us consider $n\ (\geq 3)$ variables case. 
Suppose that the following assertion (*) is true.

\medskip

(*): {\it If $f(\bbz)$ is a holomorphic weighted homogeneous polynomial with $(n-1)$ variables 
of positive degree with respect to a strictly positive weight vector $P$ 
having at least two monomials, 
then 
$f : \bc^{*(n-1)}\to \bc$ is surjective.}

\medskip

Let $m_{1} >0$ be the highest exponent of $z_{1}$. We have 
$$
f(\bbz)=p(z_{2},\dots,z_{n})z_{1}^{m_{1}}+\cdots .
$$

(1)\ Let us consider the case $f(\bbz) = p(z_{2},\dots,z_{n})z_{1}^{m_{1}}$. 
Then $p(z_2,\dots,z_{n})$ has at least $2$ monomials, and it is 
a holomorphic weighted homogeneous polynomial with $(n-1)$ variables 
of degree $d - m_{1}p_{1} > 0$ with respect to the weight vector $^t(p_{2},\dots,p_{n})$. 
By the induction hypothesis, 
$p:\bc^{*(n-1)}\to \bc$ is surjective. 
Since $f(1,z_{2},\dots,z_{n})=p(z_{2},\dots,z_{n})$, we see that $f:\bc^{*n} \to \bc$ is surjective. 

(2)\ Let us consider the case $f(\bbz)=p(z_{2},\dots,z_{n})z_{1}^{m_{1}}+\cdots+q(z_{2},\dots,z_{n})z_{1}^{n_{1}},\ q\neq 0,\ m_{1}>n_{1}\geq 0$. 
If both $p$ and $q$ are monomials, then we have $p(1,\dots,1)q(1,\dots,1) \neq 0$. 
If $p$ or $q$ has at least $2$ monomials, then 
$p(z_{2},\dots,z_{n})q(z_{2},\dots,z_{n})$ also has at least $2$ monomials, 
and 
it is a holomorphic weighted homogeneous polynomial with $(n-1)$ variables of degree $\deg pq >0$
\footnote{
Let $p(\bbz)=\sum_{\nu}b_{\nu}\bbz^{\nu},\ q(\bbz)=\sum_{\mu}c_{\mu}\bbz^{\mu}$ be 
weighted homogeneous polynomials with respect to the same weight vector $\ba$. 
Then, 
$p(\bbz)q(\bbz)=\sum_{\nu,\mu}b_{\nu}c_{\mu}\bbz^{\nu+\mu}$ is also 
a weighted homogeneous polynomials of degree $\deg p+\deg q$ with respect to $\ba$. 
Actually, 
if $\nu+\mu = \lambda$ for a fixed exponent $\lambda = (\lambda_1, \dots , \lambda_n) \in \bz_+^n$ and $b_{\nu}c_{\mu} \neq 0$, then we have 
$b_{\nu}\neq 0,\ c_{\mu} \neq 0$, and hence, 
$\sum_{j=1}^{n}a_{j}\lambda_j = 
\sum_{j=1}^{n}a_{j}(\nu_{j}+\mu_{j}) = 
\sum_{j=1}^{n}a_{j}\nu_{j} + \sum_{j=1}^{n}a_{j}\mu_{j} = \deg p+\deg q$. 
}. 
Hence, by the induction hypothesis, there exists 
$(\alpha_{2},\dots,\alpha_{n}) \in \bc^{*(n-1)}$ such that 
$p(\alpha_{2},\dots,\alpha_{n})q(\alpha_{2},\dots,\alpha_{n}) \neq 0$. 
Thus, in both cases we see that there exists $(\alpha_{2},\dots,\alpha_{n}) \in \bc^{*(n-1)}$ such that 
$p(\alpha_{2},\dots,\alpha_{n}) \neq 0$ and $q(\alpha_{2},\dots,\alpha_{n}) \neq 0$. 
Since 
$f(z_{1},\alpha_{2},\dots,\alpha_{n})=z_{1}^{n_{1}}(p(\alpha_{2},\dots,\alpha_{n})z_{1}^{m_{1}-n_{1}}+\cdots+q(\alpha_{2},\dots,\alpha_{n}))$, 
we see that 
$f:\bc^{*n} \to \bc$ is surjective by the same argument as in Step 1. %%% $2$変数の場合と同様の論法

\medskip

{\bf (iii)}\ 
We set $d_{r} := \rdeg_P f,\ d_{p} := \pdeg_Q f$. 
We will show that there are no mixed critical points of $f$ on $\bc^{*n}$. 
Note that a $0$-dimensional face (vertex) of $\Gamma_+(f)$ possibly corresponds to plural terms (monomials) of $f$.  %%%

Since $f$ is a mixed weighted homogeneous polynomial (Definition \ref{mixed-whp}) with respect to 
a radial weight vector $P$ and a polar weight vector $Q$, we have the Euler equalities:
$$
d_{r}\cdot f(\bbz,\overline{\bbz})=\sum_{j=1}^{n}\left(p_jz_{j}\frac{\partial f}{\partial z_{j}}(\bbz,\overline{\bbz})+p_j\overline{z_{j}}\frac{\partial f}{\partial
\overline{z_{j}}}(\bbz,\overline{\bbz})\right),\ \ 
d_{p}\cdot f(\bbz,\overline{\bbz})=\sum_{j=1}^{n}\left(q_jz_{j}\frac{\partial f}{\partial z_{j}}(\bbz,\overline{\bbz})-q_j\overline{z_{j}}\frac{\partial f}{\partial
\overline{z_{j}}}(\bbz,\overline{\bbz})\right) .
$$
Suppose that $\bbz^{0} \in \bc^{*n}$ is a mixed critical point of $f$. 
Then, by Proposition 1 of Oka \cite{Oka2008}, 
there exists $\alpha \in \bc$ with $|\alpha|=1$ such that
$$
\left(\overline{\frac{\partial f}{\partial
z_{1}}}(\bbz^{0},\overline{\bbz^{0}}),\cdots,\overline{\frac{\partial f}{\partial
z_{n}}}(\bbz^{0},\overline{\bbz^{0}})\right)=\alpha \left(\frac{\partial f}{\partial
\overline{z_{1}}}(\bbz^{0},\overline{\bbz^{0}}),\cdots,\frac{\partial f}{\partial
\overline{z_{n}}}(\bbz^{0},\overline{\bbz^{0}})\right) .
$$
We have $f(\bbz^{0},\overline{\bbz^{0}}) \neq 0$ since $f$ is Newton non-degenerate. 
Then we have 
\begin{eqnarray*}
d_{r}\cdot f(\bbz^{0},\overline{\bbz^{0}}) &=& \sum_{j=1}^{n}\left(p_jz_{j}\frac{\partial
f}{\partial z_{j}}(\bbz^{0},\overline{\bbz^{0}})+p_j\overline{z_{j}}\cdot \overline{\alpha}
\overline{\frac{\partial f}{\partial z_{j}}}(\bbz^{0},\overline{\bbz^{0}})\right) \\
&=& \sum_{j=1}^{n}p_j\left(z_{j}\frac{\partial f}{\partial
z_{j}}(\bbz^{0},\overline{\bbz^{0}})+\overline{\alpha}\cdot \overline{z_{j}\frac{\partial
f}{\partial z_{j}}}(\bbz^{0},\overline{\bbz^{0}})\right) \ \text{and}
\end{eqnarray*}
\begin{eqnarray*}
d_{p}\cdot f(\bbz^{0},\overline{\bbz^{0}}) &=& \sum_{j=1}^{n}\left(q_jz_{j}\frac{\partial
f}{\partial z_{j}}(\bbz^{0},\overline{\bbz^{0}})-q_j\overline{z_{j}}\cdot \overline{\alpha}
\overline{\frac{\partial f}{\partial z_{j}}}(\bbz^{0},\overline{\bbz^{0}})\right) \\
&=& \sum_{j=1}^{n}q_j\left(z_{j}\frac{\partial f}{\partial
z_{j}}(\bbz^{0},\overline{\bbz^{0}})-\overline{\alpha} \cdot \overline{z_{j}\frac{\partial
f}{\partial z_{j}}}(\bbz^{0},\overline{\bbz^{0}})\right) .
\end{eqnarray*}
We now set 
$c:=\sum_{j=1}^{n}p_jz_{j}\frac{\partial f}{\partial z_{j}}(\bbz^{0},\overline{\bbz^{0}}),\ \ 
c':=\sum_{j=1}^{n}q_jz_{j}\frac{\partial f}{\partial z_{j}}(\bbz^{0},\overline{\bbz^{0}})
$. 
Then we have 
$d_{r}\cdot f(\bbz^{0},\overline{\bbz^{0}}) = c+\overline{\alpha c}$ and 
$d_{p}\cdot f(\bbz^{0},\overline{\bbz^{0}}) = c'-\overline{\alpha c'}$. 
We set $\alpha=e^{i\theta},\ c=\rho e^{i\tau},\ c'=\rho' e^{i\tau'}$ and we have 
$c+\overline{\alpha c} = 2\rho \cos \frac{2\tau+\theta}{2}\left(\cos \left(-\frac{\theta}{2}\right)+i\sin \left(-\frac{\theta}{2}\right) \right)$ and 
$c'-\overline{\alpha c'} = 2\rho' \sin \frac{2\tau'+\theta}{2}\left(\cos \left(\frac{\pi}{2}-\frac{\theta}{2}\right)+i\sin \left(\frac{\pi}{2}-\frac{\theta}{2}\right)\right)$. 
However, 
since $d_{r},\ d_{p}$ are non-zero integers, 
the difference of the arguments of the complex numbers 
$d_{r}\cdot f(\bbz^{0},\overline{\bbz^{0}})$ and $d_{p}\cdot f(\bbz^{0},\overline{\bbz^{0}})$
is written as $n \pi \ (n \in \bz)$. 
This is a contradiction. 

\medskip

{\bf (iv)}\ Suppose that $f(\bbz)=0$ for every $\bbz \in \bc^{*n}$. 
Then $\bbz$ is a mixed critical point of $f$ and $f(\bbz)=0$. 
This contradicts the Newton non-degeneracy of $f$. 
Hence, 
there exists $\bba=(\alpha_{1},\cdots.\alpha_{n})\in \bc^{*n}$ such that $f(\bba) \neq 0$. 
Since $f$ is mixed  weighted homogeneous and $d_{r}>0,\ d_{p} \neq 0$, by Lemma \ref{mixed-action}, we have 
$$
f(t \circ (e^{i\theta} \circ \bba))=t^{d_{r}}e^{i d_{p}\theta}f(\bba)\ \ (\forall t>0,\ \forall \theta \in \br) ,
$$
and hence, $f(\bc^{*n}) \supset \bc^{*}$. 
With the assumption $f^{-1}(0) \cap {\bc^*}^n \neq \emptyset$, 
we see that 
$f : {\bc^*}^n \to \bc$ is surjective. 
This completes the proof of Proposition \ref{remark4}. \hfill \QED

%%%%%%%%%%%%%%%%%%%%%%%%%%%%%%%%%%%%%%%%%%%%%%%%%%%%%%
\medskip

Here we note that the non-zero polar degree condition ($d_{p} \neq 0$) is a key 
in both proofs of 
{\bf (iii)} (the non-existence of mixed critical points of $f$ on $\bc^{*n}$) and 
{\bf (iv)} (the surjectivity of $f : {\bc^*}^n \to \bc$) 
of Proposition \ref{remark4}. 
%%%%%%%%%%%%%%%%%%%%%%%%%%%%%%%%%%%%%%%%%

\bigskip

\begin{lemma}[cf. \cite{Oka2015},\ \S 2, p.174]\label{convenient-function-criterion}
A mixed function germ $(f,\zero)$ is convenient if and only if 
for every fixed $i \ (=1, \dots , n)$, 
there exists a pair of non-negative integers $(\nu_i,\ \mu_i)$ with $\nu_i + \mu_i >0$ 
such that 
$f$ has the term $c_{\nu, \mu}z_i^{\nu_i}\bar{z}_i^{\mu_i} \ (c_{\nu, \mu}\neq 0)$, 
where $\nu := (0, \dots , 0, \nu_i , 0, \dots , 0),\ \mu := (0, \dots , 0, \mu_i , 0, \dots , 0)$. 
\end{lemma}
\begin{proof}
This is well-known. Here we give a proof. 
Suppose that $(f,\zero)$ is convenient. 
If for some $i \ (=1, \dots , n)$, 
there are no pair of non-negative integers $(\nu_i,\ \mu_i)$ with $\nu_i + \mu_i >0$ 
such that $c_{\nu, \mu}z_i^{\nu_i}\bar{z}_i^{\mu_i} \ (c_{\nu, \mu}\neq 0)$, 
where $\nu := (0, \dots , 0, \nu_i , 0, \dots , 0),\ \mu := (0, \dots , 0, \mu_i , 0, \dots , 0)$, 
is a term of $f$, %%%% z_i と \bar{z}_i だけの積であるような項が存在しない．必ずほかの変数z_jまたは\bar{z}_j (j\neq i)が掛けられている
then we have 
$f^{ \{ i \}} \equiv 0$. 
This contradicts the convenience of $(f,\zero)$. 
Conversely, suppose that 
for every fixed $i \ (=1, \dots , n)$, 
there exists a pair of non-negative integers $(\nu_i,\ \mu_i)$ with $\nu_i + \mu_i >0$ 
such that 
$f$ has the term $c_{\nu, \mu}z_i^{\nu_i}\bar{z}_i^{\mu_i} \ (c_{\nu, \mu}\neq 0)$, 
where $\nu := (0, \dots , 0, \nu_i , 0, \dots , 0),\ \mu := (0, \dots , 0, \mu_i , 0, \dots , 0)$. 
Then 
$f$ is a mixed function 
$\sum_{\nu,\mu}c_{\nu,\mu}z_i^{\nu_i}\bar{z}_i^{\mu_i}$ with $1$ variable $z_i$ 
on 
$\bc^{ \{ i \}}\ \ (z_i \text{-axis})$. %%
By the uniqueness of the Talor expansion of a mixed function, %%
$f^{ \{ i \}} \not\equiv 0$. %%
Hence, we have $f^{I} \not\equiv 0$ for every $I \subset \{ 1,2,\cdots,n \}$ with $I \neq \emptyset$. 
Namely, $(f,\zero)$ is convenient. 
\end{proof}

%%%%%%%%%%%%%%%%%%%%%%%%%%%%%%%%%%%%%%%

\end{document}